\documentstyle{amltd2004}
\begin{document}
\annalsline{156}{2002}
\received{July 16, 1999}
\startingpage{103}
\def\bye{\end{document}}
 \font\tenrm=cmr10
\def\ritem#1{\item[{\rm #1}]}

\catcode`\@=11
\font\twelvemsb=msbm10 scaled 1100
\font\tenmsb=msbm10
\font\ninemsb=msbm10 scaled 800
\newfam\msbfam
\textfont\msbfam=\twelvemsb  \scriptfont\msbfam=\ninemsb
  \scriptscriptfont\msbfam=\ninemsb
\def\msb@{\hexnumber@\msbfam}
\def\Bbb{\relax\ifmmode\let\next\Bbb@\else
 \def\next{\errmessage{Use \string\Bbb\space only in math
mode}}\fi\next}
\def\Bbb@#1{{\Bbb@@{#1}}}
\def\Bbb@@#1{\fam\msbfam#1}
\catcode`\@=12

 \catcode`\@=11
\font\twelveeuf=eufm10 scaled 1100
\font\teneuf=eufm10
\font\nineeuf=eufm7 scaled 1100
\newfam\euffam
\textfont\euffam=\twelveeuf  \scriptfont\euffam=\teneuf
  \scriptscriptfont\euffam=\nineeuf
\def\euf@{\hexnumber@\euffam}
\def\frak{\relax\ifmmode\let\next\frak@\else
 \def\next{\errmessage{Use \string\frak\space only in math
mode}}\fi\next}
\def\frak@#1{{\frak@@{#1}}}
\def\frak@@#1{\fam\euffam#1}
\catcode`\@=12

\newcommand{\del}{\partial}
\newcommand{\bx}{\sharp}
\newcommand{\mathbb}{\bf}

\font\emi= cmmi10 scaled 1700 

\title{Rational maps are {\emi d}-adic Bernoulli}
\def\titleheadline#1{\def\one{#1}\ifx\one\empty\else
\gdef\thetitle{{\frenchspacing%
\let\\ \relax
{#1}}}\fi}
\newif\ifshort
\def\shortname#1{\global\shorttrue\xdef
\theauthors{{\eightsc\uppercase{#1}}}}
\let\shorttitle\titleheadline
\shorttitle{ \eightsc\uppercase{Rational maps are } 
{\eightpoint \it d}\eightsc\uppercase{-adic Bernoulli} }

 \twoauthors{Deborah Heicklen}{Christopher Hoffman}

\institutions{3200 Zanker Rd., MS X-75, San Jose, CA 95134
\\
{\eightpoint {\it E-mail address\/}: heicklen@math.berkeley.edu
}\\ \vglue6pt
Universitiy of Washington, Seattle, WA\\
{\eightpoint {\it E-mail address\/}: hoffman@math.washington.edu}}

\centerline{\bf Abstract}

\vglue6pt
 Freire, Lopes and Ma\~{n}\'{e} proved that for any
rational map $f$ there exists
a natural invariant measure $\mu_f$ \cite{FLM}.
Ma\~{n}\'{e} showed 
there exists an $n>0$ such that $(f^n, \mu_f)$ is
measurably conjugate to the one-sided $d^n$-shift, with Bernoulli measure
$(\frac 1{d^n},\ldots ,\frac 1{d^n})$ \cite{M2}. In this paper we show that
 $(f,\mu_f)$ is
conjugate to the one-sided Bernoulli $d$-shift. This verifies a
conjecture of Freire, Lopes and Ma\~{n}\'{e} \cite{FLM} and Lyubich
\cite{L}.  
 
\vglue-24pt
\section{Introduction}

Let $f(z)=P(z)/Q(z)$ 
be a rational map of the Riemann sphere, ${\bar C}$, of degree
$d\geq 2$. In \cite{FLM}, Freire, Lopes and Ma\~{n}\'{e} proved the existence
of a natural invariant measure $\mu_f$ for the map $f$. Namely, the
measure is the asymptotic distribution of preimages of any point $z$, except at most two exceptional points. 
Furthermore, it is the unique measure of maximal
entropy. The uniqueness was shown by both Lyubich \cite{L} and Ma\~{n}\'{e}
\cite{M1}. These properties 
and others 
are explicitly laid out in
Section \ref{flm}. 
 Freire, Lopes and Ma\~{n}\'{e} \cite{FLM} and Lyubich \cite{L} all
conjectured that the system $(f,\mu_f)$ is
conjugate to the one-sided Bernoulli $d$-shift. In this paper we give an
affirmative answer to this conjecture.

Let $X=\{0,\ldots ,d-1\}^{\mathbb N}$, ${\cal B}$ be the Borel
$\sigma$-algebra and $\sigma$ be the shift. In other words,
$\sigma(x)_n=x_{n+1}$. Let $\mu$ be product measure, where the weight
on each of the $d$ symbols is uniform, namely $\frac 1d$. This system
is called the {\it  one\/{\rm -}\/sided Bernoulli $d$\/{\rm -}\/shift}. 
An endomorphism $f$ on a measure space $(Y, {\cal C},\nu)$ is {\it $d$\/{\rm -}\/adic
Bernoulli} if $f$ is
measurably conjugate 
to the one-sided Bernoulli
$d$-shift.  
This is equivalent to the existence of
a partition $P$ of $Y$ into $d$ sets of
equal measure such that $\{f^{-n}P\}_{n\geq 0}$ are mutually independent and
$P$ generates ${\cal C}$, i.e., $\bigvee_{n\geq 0} f^{-n}P={\cal C}$. 
An example of a rational map that is $d$-adic
Bernoulli is $f(z)=z^d$.  The invariant measure $\mu_f$ is normalized
arclength on the unit circle, $\{z\ |\ |z|=1\}$. 
More generally, an endomorphism is {\it one-sided Bernoulli} if it is
conjugate to some one-sided Bernoulli shift, where the weights on the
symbols may vary.

Any measure-preserving endomorphism (noninvertible map) has a natural measure theoretic two-sided (invertible)
extension.  The two-sided extension of the one-sided Bernoulli $d$ shift is the 
{\it Bernoulli $d$\/{\rm -}\/shift}, the shift automorphism on
${\bar X}=\{0,\ldots ,d-1\}^{\mathbb Z}$ with uniform product measure on
${\bar X}$.
An automorphism (invertible map) is
{\it Bernoulli} if it is conjugate to a Bernoulli shift. 
If an endomorphism is one-sided Bernoulli then its two-sided extension is
Bernoulli.  The converse is not true.  For example, there exist Markov
endomorphisms which are not one-sided Bernoulli but whose two-sided
extension are Bernoulli \cite{AMT}, \cite{FO}.

Historically, determining whether an endomorphism is one-sided Bernoulli
is a much harder question than whether an automorphism is
Bernoulli. This is because there is a reasonable condition for
determining whether an automorphism is Bernoulli called very weak
Bernoulli (VWB) \cite{O}.  Ornstein and Weiss proved that an
automorphism is Bernoulli if and only if it is very weak Bernoulli
\cite{O}, \cite{OW}. 
Many natural systems have been shown to be either Bernoulli
or not Bernoulli by checking this condition.  Examples of Bernoulli
systems are toral automorphisms and geodesic flows on spaces of
constant negative
curvature \cite {K}, \cite{OW1}.  In the other direction, 
Kalikow showed that the $[T,T^{-1}]$
system is not Bernoulli \cite{Kal}.

Recently, Hoffman and Rudolph  developed a condition
analogous to VWB called tree very weak Bernoulli (tree vwB) \cite{HR}. 
In this paper
they showed that tree vwB is equivalent to one-sided Bernoulli. We define
and explain this condition in Section \ref{wroo}. With this condition
there is no need to construct a conjugacy  
or to find a $d$ set independent generating partition.

In \cite{M2} Ma\~{n}\'{e} proved 
that there exists an $n>0$ such that the system
$(f^n,\mu_f)$ is conjugate to the one-sided Bernoulli
$d^n$-shift. (Throughout this paper, $f^n$
denotes iteration.)
Since $(f^n, \mu_f)$ is one-sided Bernoulli, its two-sided extension
is isomorphic to a Bernoulli automorphism.  By Ornstein's theory, the
two-sided extension of $(f, \mu_f)$ is isomorphic to a Bernoulli
shift \cite{O2}.
It is important to note that, in contrast to the case of automorphisms,
the fact that $f^n$ is one-sided Bernoulli does not imply that $f$
itself is one-sided Bernoulli \cite{H2}.

One possible way to prove that $f$ is $d$-adic Bernoulli is as
follows. Since $f$ has finitely many critical values (the image under $f$ of 
points with $f'=0$), it is possible
to connect them with an arc $\gamma$ such that $\mu_f(\gamma)=0$. Now
$f^{-1}(\gamma^c)$ is an open set with $d$ connected components. Let
$P$ be the partition whose sets are these connected components. It can
be shown that $\{f^{-i}P\}_{i\geq0}$ is independent and for every
$p\in P$, $\mu_f(p)=\frac 1d$. However it is not clear that $P$
generates. 
This is the general approach that Ma\~{n}\'{e} used in \cite{M2}.



There are two main tools we use in the proof that $(f,\mu_f)$ is $d$-adic
Bernoulli. The first is the tree vwB condition described in Section
\ref{wroo}.  The second is the construction of the unique measure of
maximal entropy
\cite{FLM}. We show that for any two points $z,z'$ in the support of this 
measure, the set $\{f^{-n}(z)\}$ can be matched in a one-to-one manner 
with the points in $\{f^{-n}(z')\}$ in such a way that the matching 
preserves the underlying tree structure of $\{f^{-n}(z)\}$ and so that most paired points are close together.
We then use this matching to show that the endomorphism is tree vwB.

\section{Rational endomorphisms}
\label{flm}

In this section we present some basic facts about rational maps.
For more information about this see \cite{C}.  Then we
describe the unique invariant measure of maximal entropy and 
some of its properties.

Let $f(z)$ be a rational endomorphism of ${\bar C}$. Then 
$f(z)=P(z)/Q(z)$ where $P$ and $Q$ are relatively prime
polynomials. Define the {\it degree of $f$} to be the maximum degree
of $P$ and $Q$. Throughout we assume the degree of $f$ 
is at least 2.  We use
${\rm dist}(x,y)$ to represent the distance between $x$ and $y$ on the
Riemann sphere.  Thus ${\rm dist}(x,y)$ is bounded by 1.

The {\it Julia set}, denoted by $J(f)$ is the set of all
points $z\in{\bar C}$ such that for every neighborhood $U$ of $z$,
$\{f^n|_U\}$ is not a normal family. That is, no subsequence of 
this family of functions converges uniformly on compact subsets. 
$J(f)$ is a nonempty compact invariant set. Moreover, it is the closure of the repelling
periodic orbits. 
The simplest example of a Julia set is for the map $f(z)=z^2$ when  $J(f)$ is the unit circle.
For the map $f(z)=1-2/z^2$  the Julia set is $\bar C$.
However for other rational functions the Julia set can be very complex.
For the map $f(z)=z^2+1$, $J(f)$ is a totally disconnected set.  In this case it is also fractal and 
conformally self-similar.  There are also examples of rational functions for which the Julia is connected but not
locally connected, as well as functions for which the Julia set is neither connected nor totally disconnected.
  
Fatou and Julia introduced a set Exc$(f)\subset {\bar C}$.  
This set is the maximal finite set which is invariant under $f$ and $f^{-1}$. 
They proved it contains at most two points. 
The measure $\mu_f$ is
the weak star limit of measures uniformly supported on $\{f^{-n}(a)\}$
for any $a\notin$ Exc$(f)$. 
Brolin first introduced this measure for the case when $f$ is a polynomial \cite{B}.  Lyubich \cite{L1} and Freire, Lopes and Ma\~{n}\'{e} \cite{FLM} generalized
this to the case of rational functions.
This limit
does not depend on the choice of $a$. More precisely, define the {\it
$n$\/{\rm -}\/preimages of $a$} to be the set $\{f^{-n}(a)\}$. (The cardinality
of $\{f^{-n}(a)\}$ is $d^n$ as elements are listed possibly multiple
times according to their multiplicities.  We will use this convention throughout the paper.) Order the $n$-preimages of
$a$ as $z_i^n(a)$, $i=1,\ldots ,d^n$. Define
$$\mu_n(a)=\frac 1{d^n}\sum \delta_{z_i^n(a)}$$ where
$\delta_{z_i^n(a)}$ is Dirac measure supported on $z_i^n(a)$. 
The space of invariant probability measures on ${\bar
C}$ is endowed with the weak star topology. 

\proclaimtitle{\cite{FLM}}
\proclaim{Theorem}  There exists an $f$ invariant probability
measure $\mu_f$ satisfying the following properties.
\begin{itemize}

\ritem{1.} $\lim \mu_n(a)=\mu_f$ for all $a\notin$ {\rm Exc}$(f)$. Moreover
this convergence is uniform as $a$ varies over a compact subset of
{\rm Exc}$(f)^c$.

\ritem{2.} Support $\mu_f=J(f)$.
\ritem{3.} $f$ is exact.
\ritem{4.} For all Borel sets $A$ such that $f|_A$ is injective{\rm ,} $\frac 1d
\mu_f(fA)=\mu_f(A)$. 
\ritem{5.} $h_{\mu_f}(f)=\log d$.
\end{itemize}
\endproclaim 

Properties 3--5 are conjugacy invariants. They are necessary conditions
in order for $(f, \mu_f)$ to be measure theoretically isomorphic to the one-sided Bernoulli
$d$-shift, which also satisfies these properties. 
The endomorphism $f$ is exact if $\mathbold{\cap}_{n\geq 0}f^{-n}({\cal A})$
is the trivial $\sigma$-algebra, where ${\cal A}$ is the Borel
$\sigma$-algebra on ${\bar C}$.  This property implies the endomorphism is mixing but is
much stronger. 
The last property listed implies that $\mu_f$ is a
measure of maximal entropy. This follows from Gromov's result that
the topological entropy of $f$ is $\log d$ \cite{G}. 
Lyubich \cite{L}, \cite{L2} and Ma\~{n}\'{e} \cite{M1} showed 
that $\mu_f$ is the unique measure of maximal entropy. 

\section {Tree very weak Bernoulli}
\label{wroo}

In this section we outline the ergodic theory necessary for this paper.
Let $T$ be a measure-preserving endomorphism on a probability space
$(Y,{\cal C}, \nu)$. It causes no loss of generality to assume that $Y$ is a compact metric space with metric $D$.
Also, $T$ is $d$-\/{\it adic} if almost every atom of $T^{-1}{\cal C}$
consists of $d$ points in ${\cal C}$, and the conditional measure of these 
points given the atom is $1/d$  almost everywhere. For completeness we include the proof of the following well known fact.

\proclaim{Lemma}  If $f$ is a rational map of degree $d$ then $(f,\mu_f)$ is a $d$\/{\rm -}\/adic system.
\endproclaim 

\demo{Proof} In order to show $(f,\mu_f)$ is $d$-adic, it
suffices to show that there exists a $d$ set partition $P$ satisfying
the following properties. Let ${\cal A}$ be the Borel $\sigma$-algebra of $\bar{C}$.
For all $p\in P$, $\mu_f(p)=\frac 1d$ and
$P$ and $f^{-1}{\cal A}$ are independent. Consider the $d$ set
partition $P$ constructed in the introduction. For any $A\in {\cal
A}$ and  $p\in P$, it suffices to show
\begin{eqnarray*}
\mu_f(f^{-1}A\cap p) &=& \mu_f(p)\mu_f(f^{-1}A)\\ &=& \frac 1d \mu_f(A).
\end{eqnarray*}

This follows from the fact that $f$ is injective on $f^{-1}A\cap p$
and that\break $f(f^{-1}A\cap p)=A$. Hence, $\mu_f(f(f^{-1}A\cap
p))=d(\mu_f(f^{-1}A\cap p))$ and the desired equality follows. \enddemo

For any $d$-adic system, define an {\it $N$\/{\rm -}\/tree of a point $z$} to be a
$d$-ary tree with $N+1$ levels, 0 through $N$, whose $d^n$ vertices at the $n$-th level are
identified with the $d^n$ elements of $T^{-n}(z)$. Furthermore, the
vertices are labeled in such a way that if $x$ and $y$ have the same
parent vertex, then $f(x)=f(y)$. 
We define an {\it automorphism of the $N$\/{\rm -}\/tree} of $z$ to be a map $A$ 
from the vertices of the $N$-tree of $z$ to themselves that preserves the tree structure. In
particular, if $x=f(y)$, then $A(x)=f(Ay)$.  
The labeling of the vertices is not unique because applying an automorphism of the
$N$-tree of $z$ results in a different labeled tree which satisfies
the definition.  This ambiguity will not make a difference.
Define a {\it tree of a point
$z$} to be an infinite $d$-ary tree, whose vertices are identified with
the preimages of $z$ in the same way as before. 

Given two points $z$ and $w$ we now define a metric between the $N$ trees of $z$ and $w$.
For each $n$ choose an ordering of $T^{-n}(z)$ and $T^{-n}(w)$. Call them $z^n_i$ and $w^n_i,$ $i=1,\ldots ,d^n$. 
Set

$$t_N^D(z,w)= \min_A \frac 1{N} \sum_{n=1}^N \frac 1{d^n}\sum_{i=1}^{d^n}
	{D(z^n_i, A(w^n_i))}$$ 
where $A$ ranges over all
automorphisms of the $N$-tree of $w$.  
Each automorphism $A$ generates a pairing of the vertices of the $N$-trees of $z$ and $w$.  The quantity that we are minimizing 
is a weighted average of the distance between paired vertices.  Thus this metric is a minimum distance between two sets of preimages,
subject to the constraint that the pairing of the two sets respects the tree structures of the two sets.

Another way to view this quantity is to notice that for any $z,w,$ and $A$ 
$$\frac 1{N} \sum_{n=1}^N \frac 1{d^n}\sum_{i=1}^{d^n} {D(z^n_i, A(w^n_i))}
	=\frac 1{d^N} \sum_{i=1}^{d^N} \frac 1{N}\sum_{j=0}^{N-1}
	{D(T^j(z^N_i), T^j(A(w^N_i)))}.$$ 
Each automorphism $A$ generates a pairing of $T^{-N}(z)$ and $T^{-N}(w)$.  The first summand is the average distance
between images of $z_i^N$ and $A(w_i^N)$ under $T^0,T^1,\ldots ,T^{N-1}$. The outer
\pagebreak sum averages this quantity over $T^{-N}(z)$.
Viewed in this manner the definition of tree vwB is close to the definition of very
weak Bernoulli for two-sided extensions of $d$-adic endomorphisms.  
See \cite{HR} for a discussion of this definition difference.

\specialnumber{3.1}\numbereddemo{Definition} 
A $d$-adic endomorphism $T$ acting on $(Y,{\cal C},\nu)$  is {\it tree very weak
Bernoulli} (tree vwB) if  for all $\varepsilon>0$,  there exist $N$ and a set $G$
such that
\begin{itemize}
\item[1.] $\nu(G)>1-\varepsilon$.
\item[2.] For all $z,w\in G$, $t_N^D(z,w)<\varepsilon$.
\end{itemize}
\enddemo

The main result (Theorem 5.5) of \cite{HR} is the following.

\specialnumber{3.1}\proclaim{Theorem} \label{hr} Let $T$ be a 
$d$\/{\rm -}\/adic endomorphism acting on $(Y,{\cal C},\nu)$ which is tree {\rm vwB.}
 Then $T$ acting on $(Y,{\cal C},\nu)$
is $d$\/{\rm -}\/adic Bernoulli.
\endproclaim
 
\section {Rational maps are  Bernoulli}

In this section we verify the tree vwB condition for all rational maps with the dist metric.
The main tool is Lemma \ref{fundlem}.

\numbereddemo{Definition} Two points $z,w\in{\bar C}$ are $(N,\varepsilon)$ tree
related if there exists  an invertible map
$$\phi: \mbox{vertices of the tree of $z$} \to \mbox{vertices of the tree of $w$}$$ and for each
$n$ there exists
an ordering of $f^{-n}(z)$, $\{z^n_i\}_{i=1}^{i=d^n}$, such that
\begin{itemize}
\item[1.] $\phi(f(z^n_i))=f(\phi(z^n_i))$ for all $n$ and $i$.
\item[2.] ${\rm dist}(z^n_i,\phi(z^n_i))<\varepsilon$, 
	for all $n\geq N$ and  $i,\ 1\leq i\leq(1-\varepsilon)d^n$.
\end{itemize}
\enddemo

 The main idea of this section is to show that any two points in
the Julia set are $(N,\varepsilon)$ tree related. This will imply the tree
vwB condition.
The next lemma shows the relationship between $(N',\varepsilon)$ tree related
and the $t^{{\rm dist}}_N$ metric.
\advance\theoremcount by -1
\proclaim{Lemma} \label{new}
If $z$ and $w$ are $(N',\varepsilon)$ tree related and $N>N'/\varepsilon$ then 
$$t_{N}^{{\rm dist}}(z,w)<3 \varepsilon.$$
\endproclaim 

\demo{Proof}
The map $\phi$ generates an automorphism of the $N$ tree of $w$ by 
$$A(w^{n}_{i})=\phi(z^{n}_{i}).$$
Thus if $n>N'$ and $i \leq (1-\varepsilon)d^n$ then 
$${\rm dist}(z^n_i,A(w^n_i)) < \varepsilon.$$
Thus
\begin{eqnarray*}
&&\kern-3.5in 
\frac 1{N}   \sum_{n=1}^{N} \frac1{d^n} \sum_{i=1}^{d^n} {{\rm dist}(z^n_i, A(w^n_i))}
\end{eqnarray*}
\begin{eqnarray*}
& \le & \frac 1{N}\left( \sum_{n=1}^{N'} \frac1{d^n} \sum_{i=1}^{d^n} {{\rm dist}(z^n_i, A(w^n_i))}
					+\sum_{n=N'+1}^{N} \frac1{d^n} \sum_{i=1}^{(1-\varepsilon)d^n} {{\rm dist}(z^n_i, A(w^n_i))} \right.\\
&	&				\left.+\sum_{n=N'+1}^{N} \frac1{d^n} \sum_{i=(1-\varepsilon)d^n+1}^{d^n} {{\rm dist}(z^n_i, A(w^n_i))}
				\right)\\
& \le & \frac 1{N}\left( \sum_{n=1}^{N'} \frac1{d^n} \sum_{i=1}^{d^n} 1
					+\sum_{n=N'+1}^{N} \frac1{d^n} \sum_{i=1}^{d^n} \varepsilon
					+\sum_{n=N'+1}^{N} \frac1{d^n} \sum_{i=(1-\varepsilon)d^n+1}^{d^n} 1
				\right)\\
& \le & \frac 1{N}(N'+\varepsilon N + \varepsilon N)\\
& \le & 3 \varepsilon.\\
\noalign{\vskip-36pt}
\end{eqnarray*}
\enddemo
\vglue9pt

For the remainder of the paper we say the map $g: U \to V$ is a $k$-to-one map if $g$ is a branched covering of finite
degree such that, for each $v \in V$, $k$ is equal to the sum over $g^{-1}(v) \cap U$ of the multiplicity of
the solution of $g(u)=v$.
 
\numbereddemo{Definition} A topological disk $U$ is $(N,\varepsilon)$ tree adapted
if there exist  topological disks, $S^0=U$,  $\{S^n_i\}$, $1\leq i\leq
l_n$, and integers $k_i^n$, for every $n\geq 1$ such that
\begin{itemize}

\item[1.] $\sum_{i=1}^{l_n} k^n_i=d^n$ for all $n\leq N$.

\item[2.] $\sum_{i=1}^{l_n} k^n_i\geq(1-\varepsilon)d^n$ for all $n\geq1$.

\item[3.] For all $n\leq N$ for all $j$, $f^{n}|_{S^n_j}$ is a $k_j^n$-to-one
map onto $S^0$. Moreover,  for all $n\leq N$ for all $j$, there exists $i$
such that $f|_{S^n_j}$ is a finite to one map of $S^{n-1}_i$.

\item[4.] For all $n>N$ for all $j$ there exists $i $ such that $f$ maps $S^n_j$
homeomorphically onto $S^{n-1}_i$.

\item[5.] $\lim_{n\to\infty}\sup_i$ diam $S^n_i=0$.
\end{itemize}
\enddemo

The following lemma ties together the concepts of tree adpated and tree related. 

\advance\theoremcount by -1
\proclaim{Lemma} Given an $(N,\varepsilon)$ tree adapted set $U${\rm ,} there 
exists an $M$ such that for any two points $z,w\in U${\rm ,} the points $z$ and
$w$ are $(M,\varepsilon)$ tree related.
\endproclaim 

\demo{Proof} Fix $z,w\in U$, an $(N,\varepsilon)$ tree adapted
set. There exists a collection of sets, $\{S^n_i\}$,
 such that $\lim_{n\to\infty}\sup_i$ diam$S^n_i=0$. Pick $M$
such that diam$S^n_i<\varepsilon$ for all $n\geq M$ and $i \leq l_n$. Define
$\phi(z)=w$, and define
$\phi:\{f^{-1}(z)\}\to \{f^{-1}(w)\}$ in a one-to-one manner so that if $z_i^1\in 
S^1_k$, 
then $\phi(z_i^1)\in S^1_k$.
 Define $\phi$ on the
rest of the tree by induction. Assume $\phi$ has been constructed on
the first $n$ levels such that if $\phi(z^j_i)=w^j_k$ then
$\phi(f(z^j_i))=f(w^j_k)$ for all $i,k$ and for $j\leq n$. Also
assume that if $z^j_i\in S^j_k$ then
$\phi(z^j_i)\in S^j_k$ for all $k$ and  $j\leq n$.
 For $n+1$, pick $z^{n+1}_i$ which is contained in some 
$S^{n+1}_k$. Then $f(z^{n+1}_i)\in S^n_j$ for some $j$, and
$\phi(f(z^{n+1}_i))\in S^n_j$. Notice that
$$\#\{f^{-1}(f(z^{n+1}_i))\cap S^{n+1}_k\}=\#\{f^{-1}(\phi(f(z_i^{n+1})))\cap
S^{n+1}_k\}$$ 
since $f|_{S^{n+1}_k}$ is a branched covering of finite degree. Define
$\phi$ to be any one-to-one map between these sets. On the preimages that lie
inside some $S^{n+1}_i$, $\phi$ fulfills the requirements. On
preimages not lying in one of these sets, define $\phi$ to be any map
that preserves the tree structure, namely, any map that sends
$f^{-1}(f(z^{n+1}_i))$ to $f^{-1}(\phi (f(z^{n+1}_i)))$. Thus
$\phi(f(z^{n+1}_i))=f(\phi( z^{n+1}_i))$. Furthermore,  for all $n\geq M$,  
 for all $i$ such that $z_i^n\in\mathbold{\cup}_{j=1}^{l_n} S_j^n$,
$$\ {\rm dist}(z^n_i,\phi(z^n_i))<\varepsilon,$$
since $z_i^n,\phi(z_i^n)\in S^n_j$ for some $j$  and diam$S
^n_i<\varepsilon$.   The cardinality of 
the set $\{i\ |\ z^n_i\in\mathbold{\cup}_{j=1}^{l_n}S_j^n\} $ is at 
least $(1-\varepsilon)d^n$ since $\sum_{i=1}^{l_n}k_i^n\geq(1-\varepsilon)d^n$.
Thus the desired orderings exist.
\enddemo
 
The fundamental lemma of this section is the following.

\specialnumber{4.3}
\proclaim{Lemma} 
\label{fundlem}
Given $\varepsilon>0${\rm ,} $z\notin$ {\rm Exc}$(f)${\rm ,} and an arc
$\gamma$ containing $z$ such that $\gamma\setminus\{z\}$ does not
contain any critical values of $f^n$ for all $n\geq1${\rm ,} there exists an
$(N,\varepsilon)$ tree adapted set $U$ containing $\gamma$ for some $N\geq1$.
\endproclaim 

We leave the proof of this lemma until the end of the section.  The
fundamental lemma 
gives us the following lemma.

\specialnumber{4.4}\proclaim{Lemma} \label{bounded}
Given $\varepsilon>0$ there exists $N$ such that all pairs of points
$x,y\in J(f)$ are $(N,\varepsilon)$ tree related.
\endproclaim 

\demo{Proof} First we show that for any two points $x,y\in
J(f)$, there exists an $M$ such that $x,y$ are $(M,\varepsilon)$ tree
related. Given $\varepsilon>0$ pick arcs $\gamma_x$ and $\gamma_y$
satisfying the hypothesis of the fundamental lemma, 
such that $x\in \gamma_x$, $y\in\gamma_y$,  and such that
$\gamma_x\cap\gamma_y\neq\emptyset$. By the fundamental lemma, there
exist  topological disks $U_x$ and $\ U_y$ containing $\gamma_x$ and $\gamma_y$
that are $(N,\frac\varepsilon2)$ tree adapted. Pick $M$ such that any two
points in $U_x$ (and $U_y$) are $( M,\frac\varepsilon2)$ tree
related. This implies that $x,y$ are $(M,\varepsilon)$ tree related since
$U_x\cap U_y\neq\emptyset$. In particular, suppose $z\in U_x\cap U_y$,
$\phi_x$ is the map from the tree of $z$ to that of $x$ and $\phi_y$
is the map from the tree of $z$ to the tree of $y$. Then
$\phi=\phi_x\phi_y^{-1}$ is a map from the tree of $y$ to the tree of
$x$ preserving the tree structure. Furthermore, if $x_i^n\in S_j^n$, for some $j \leq l_n$, and $n\geq M$,
then ${\rm dist}(x^n_i,\phi
(x^n_i))<\varepsilon$. This implies that the
 cardinality of 
the set $\{i\ |\ {\rm dist}(x^n_i,\phi(x_i^n))<\varepsilon\} $ is at 
least $(1-\varepsilon)d^n$ since $\sum_{i=1}^{l_n}k_i^n\geq(1-\varepsilon)d^n$.

Now define ${\tilde M}:J(f)\times J(f)\to {\mathbb N}$ by the
following property. Let ${\tilde M}(z,w)$ be the minimum $N\geq 0$
such that there are neighborhoods $U_z$ and $U_w$ such that every point in
$U_z$ is $(N,\varepsilon)$ tree related to every point in $U_w$. The
previous argument shows that this is well defined. Furthermore, this
function is upper semicontinuous. Since $J(f)$ is compact ${\tilde M}$ 
is bounded. Let $N$ be an
upper bound. $N$ satisfies the required property.  \enddemo

\specialnumber{4.1}
\proclaim{Theorem} If $f$ is a rational map of degree $d\geq 2$ then the 
system $f$ acting on $(\bar C, {\cal A},\mu_f)$ is $d$\/{\rm -}\/adic Bernoulli.
\endproclaim

\demo{Proof} 
Given $\varepsilon>0$.
Set $G=J(f)$.  By Lemma \ref{bounded} pick $N'$ such that for all $z,w \in J(f)$ are $(N',\varepsilon/3)$ tree related.  
Choose $N>3N'/\varepsilon$.
By Lemma \ref{new}
$$t^{{\rm dist}}_N(z,w)< \varepsilon.$$
Thus $f$ acting on $(\bar C, {\cal A},\mu_f)$ is tree vwB. By Theorem \ref{hr}, $f$ 
acting on $(\bar C, {\cal A},\mu_f)$ is $d$-adic Bernoulli. 
\enddemo

In order to prove the fundamental lemma we need a preliminary lemma from
\cite{FLM} that we state but do not prove. For any $z\in{\bar C}$
define 
$$m_n(z)=\mbox{ the maximum multiplicity of any }
n\mbox{-preimage of }z.$$

\specialnumber{4.5}
\proclaim{Lemma} For all $z\notin$ {\rm Exc}$(f)$ there exists $N>0$ and 
$d_0,\ 1\leq d_0< d$ such that $m_n(z)\leq (d_0)^n$ for all $n \geq N$.
\endproclaim

\demo{Proof of Lemma {\rm \ref{fundlem}}} 
Pick $N_0$ such that $m_{N_0}(z)N_0
d^{-N_0}<\frac \varepsilon2$. Such an $N_0$ exists by the previous
lemma. Furthermore, pick $N_0$ such that

$$4m_{N_0}(z)d^{-N_0}\frac{d^3}{d-1}<\frac\varepsilon4$$
and
$$2m_{N_0}(z)N_0d^{-N_0}(1+\sum_{j=1}^\infty\frac j{d^j})<\frac\varepsilon4.$$


Let $m_{N_0}(z)=m$. Since the only possible 
critical value of $f^{N_0}$, contained
in $\gamma$, is $z$, it follows that the connected components
$\{\gamma^n_i\}_{n=1}^{N_0}$ of $f^{-n}(\gamma)$ are either arcs or
unions of arcs with a unique point of intersection. Therefore, each
$\gamma^n_i$ is simply connected. We can then \pagebreak take a topological disk
$U_0\supset\gamma$ so thin that for each $n\leq N_0$ there exist 
disjoint topological disks $U_i^n$ containing $\gamma^n_i$,   such
that for all $n > N_0$, for all $i$,  there exists $j$ such that
$f(U^n_i)=U^{n-1}_j$.  Furthermore, $f|_{U^{n}_i}$ is a
$k^n_i$-to-one map (counting multiplicity) onto $U_0$.

Set $\varepsilon_{N_0}=2mN_0d^{-N_0}$ and
$\varepsilon_{n+1}=\varepsilon_n+4md^{-(n+1)}d^2+2md^{-(n+1)}(n+1)$ for
$n\geq N_0$. Observe that
\begin{eqnarray*}
\varepsilon_n 
& \leq & 4md^2\sum_{j=N_0}^\infty\frac 1{d^j}
	+2m\sum_{j=N_0}^\infty \frac j{d^j}\\ 
& \leq & 4md^{-N_0}\frac{d^3}{d-1}+2mN_0d^{-N_0}
	(1+\sum_{j=1}^\infty\frac{j}{d^j})\\ 
& \leq & \frac\varepsilon2.
\end{eqnarray*} 

We claim that  for all $n\geq N_0,$ $f^{-n}(U_0)$ contains a union of
topological disks, $W^n_i,\ i=1,\ldots ,l_n$ such that $f^n(W^n_i)=U_0$ is
a $k^n_i$-to-one map, where $k^n_i,\ l_n$ are integers satisfying

\begin{itemize}
\item[1.] $1\leq k^n_i\leq m$ and
\item[2.] $\sum_{i=1}^{l_n}k^n_i\geq (1-\frac{\varepsilon_n}2)d^n$.
\end{itemize}
 Furthermore $\lambda(W^n_i)<\frac 1n$, where $\lambda$ denotes
Lebesgue measure. For $N_0$, Let $W^n_i=U^n_i$. Since $\{U^{N_0}_i\}$
are disjoint, there are at most $N_0$ of these disks such that
$\lambda(U^{N_0}_i)\geq\frac 1{N_0}$. Throw these disks away, 
so that if $1\leq i\leq l_{N_0}$, then $\lambda(U^{N_0}_i)<\frac
1{N_0}$. Now
\begin {eqnarray*}
\sum_{i=1}^{l_n} k^{N_0}_i 
& \geq & d^{N_0}-N_0m=d^{N_0}(1-mN_0d^{-N_0})\\
& \geq & d^{N_0}(1-\frac{\varepsilon_{N_0}}{2}).
\end{eqnarray*}

The proof of the claim is completed by induction. Suppose there exists
the collection of topological disks $W^n_i$ and integers $k^n_i,\
i=1,\ldots ,l_n$. Let $H$ be the set of integers $t$ between 1 and $l_n$
such that $W_i^n$ contains no critical values of 
$f$. For every $t\in H$ there is a disk that maps homeomorphically
onto $W_t^n$. Define this disk as $W^{n+1}_i$ and only keep those $i$
such that $\lambda(W^{n+1}_i)<\frac 1{n+1}$. These will be the values
$1\leq i\leq l_{n+1}$.  Then
\begin{eqnarray*}
\sum_{i=1}^{l_{n+1}}k^{n+1}_i 
 & \geq & d(\sum_{i\in H}k^n_i)-(n+1)m\\
 & \geq & d(\sum_{i=1}^{l_n}k^n_i-\sum_{i\notin H}k^n_i)-(n+1)m\\
 & \geq & d\sum_{i=1}^{l_n}k^n_i-dm(l_n-\#H)-m(n+1)\\
 & \geq & d^{n+1}(1-\frac{\varepsilon_n}2) -dm(l_n-\# H)-m(n+1).
\end{eqnarray*}
But $l_n-\# H$ is bounded by the number of critical values of $f$,
which is bounded by $2d$. Hence
\begin{eqnarray*}
  \sum_{i=1}^{l_{n+1}}k^{n+1}_i 
 & \geq & d^{n+1}(1-\frac{\varepsilon_n}2)-2d^2m-m(n+1)\\
& \geq & d^{n+1}(1-\frac 12(\varepsilon_n+4d^2md^{-(n+1)}+2m(n+1)d^{-(n+1)})\\
& \geq & d^{n+1}(1-\frac{\varepsilon_{n+1}}2).
\end{eqnarray*}

By the way, the $W^n_i$ are constructed, for $n>N_0$, for all $i$, there
exists a $j$ such that $f$ restricted to $W^n_i$ is a homeomorphism
onto $W^{n-1}_j$. Furthermore,   if we omit some $W_i^n$
because the diameter of the set is too large or because $W_i^n$
contains a critical value of $f$, then no subsequent
preimage of the set is included in the collection. We have now shown
the first four properties. In order to show the last property, we need
to shave down the set $U_0$ to a set $U$  and the sets $W^n_i$ to sets
$S^n_i$ and apply Koebe's distortion theorem. We do this exactly as in
\cite{FLM}. The proof is included for completeness.

We show that for any topological disk $U$ whose closure
 is contained in $U_0$,
$$\lim_{n\to\infty}(\sup_i\mbox {diam}(f^{-n}(U)\cap W^n_i))=0.$$

If this property is true, the lemma is proved taking $U$ containing
$\gamma$ and with closure contained in $U_0$. Then we define
 $$S^n_i=W^n_i\cap f^{-n}(U).$$
By the way the sets $W_i^n$ are  constructed,  $f^{n-N_0}|_{W^n_i}$ is
a conformal representation onto some $W^{N_0}_j$. Let
$\phi^n_i:W^{N_0}_j\to W^n_i$ be its inverse. Set $D_r=\{z\ |\ |z|\leq
r\}$. Let $\alpha_j: D_1\to W^{N_0}_j$ be a conformal representation. 
Define $\psi^n_i:D_1\to W^n_i$ as
$\psi^n_i=\phi^n_i\alpha_j$. We shall prove 
$$\lim_{n\to\infty}(\sup_i\mbox {diam}\psi_i^n(D_r))=0$$
for all $r, \ 0<r<1$. This implies the result because
$\psi^n_i(D_r)\supset f^{-(n-N_0)}(U)\cap W^n_i$ for all $n\geq N$, if
$r$ is near enough to 1. To prove the result, recall Koebe's
distortion theorem for univalent functions.  This theorem says that for 
all $0<r<1$ there exists $K(r)$ such that every univalent function
$\phi:D_1\to {\mathbb C}$ satisfies $|\phi'(a)/\phi'(b)|\leq K(r)$ for
all $a$ and $b$ in $D_r$. In particular,
$$\lambda(\phi(D_r))\geq K(r)^{-1}|\phi'(a)|\lambda(D_r)$$
for all $a\in D_r$. In our case
$$\frac 1n\geq \lambda(W^n_i)\geq \lambda(\psi^n_i(D_r))\geq K(r)^{-1}\lambda(D_r) |(\psi^n_i)'(z)|$$
for all $0<r<1,\ z\in D_r$. Then 
$$\lim_{n\to \infty}\sup_{i,z\in D_r}|(\psi^n_i)'(z)|=0$$
which implies the result.
\enddemo

 \vglue-12pt


\begin{references} 
\bibitem{AMT}
\name{J. Ashley, B. Marcus}, and \name{S. Tuncel},
The classification of one-sided Markov chains,
{\it Ergod.\ Theory Dynam.\ Systems\/}
{\bf 17}
(1997),
269--295.


\bibitem{B} \name{H. Brolin},
Invariant sets under iteration of rational functions,
{\it Ark.\ Math\/}.\
{\bf 6}
(1965),
103--144.


\bibitem{C} \name{L. Carleson} and \name{T. Gamelin},
{\it Complex Dynamics\/}, Springer-Verlag, New York, 1993.

\bibitem{FO}  \name{N. Friedman} and \name{D. Ornstein}, 
On isomorphism of weak Bernoulli transformations, 
{\it Adv.\ in Math\/}.\ 
{\bf 5}
(1970),
365--394.

\bibitem{FLM} \name{A. Freire, A. Lopes}, and \name{R. Ma\~{n}\'{e}}, 
An invariant measure for rational maps, 
{\it Bol.\ Soc.\ Brasil.\ Mat\/}.\ 
{\bf 14}
(1983),
 45--62.

\bibitem{G} \name{M. Gromov},
On the entropy of holomorphic maps,
 IHES preprint.


\bibitem{HR} \name{C. Hoffman} and \name{D. Rudolph},
Uniform endomorphisms which are isomorphic to a Bernoulli shift,
{\it Ann.\ of Math\/}.\ {\bf 156} (2002), 79--101.

\bibitem{H2} \name{C. Hoffman},
A dyadic endomorphism whose square is Bernoulli,
manuscript in progress.

\bibitem{Kal}  \name{S. Kalikow}, 
$T,T^{-1}$ transformation is not loosely Bernoulli,
{\it Ann.\  of Math\/}.\ 
{\bf 115} 
(1982), 
393--409.

\bibitem{K} \name{Y. Katznelson},
Ergodic automorphisms of $T^n$ are Bernoulli shifts,
{\it Israel J.\ Math\/}.\
{\bf 10}
(1971),
186--195.




\bibitem{L} \name{M. Lyubich},
Entropy properties of rational endomorphisms of the Riemann sphere,
{\it Ergod.\   Theory Dynam.\ Systems}
{\bf 3}
(1983),
351--385.


\bibitem{L1} \bibline,
The measure of maximal entropy of a rational endomorphism of a Riemann
sphere,
{\it Funct.\ Anal.\   Appl\/}.\
{\bf 16}
(1982),
309--311.



\bibitem{L2} \bibline,
Entropy of analytic endomorphisms of the Riemann sphere,
{\it Funct.\  Anal.\  Appl\/}.\
{\bf 15}
(1981),
300--302.

\bibitem{M1} \name{R. Ma\~{n}\'{e}}, 
On the uniqueness of the maximizing measure for rational maps,
{\it Bol.\ Soc.\ Bras.\ Mat\/}.\
{\bf 14}
(1983),
27--43.


\bibitem{M2} \bibline,
On the Bernoulli property of rational maps,
{\it Ergod.\ Theory Dynam.\ Systems}
{\bf 5}
(1985),
71--88.

\bibitem{O}  
\name{D. Ornstein}, 
Two Bernoulli shifts with infinite entropy are isomorphic,
{\it Adv.\ in Math\/}.\  
{\bf 5}
(1970),
339--348.

\bibitem{O2}           
\bibline,
Factors of Bernoulli shifts are Bernoulli shifts,
{\it Adv.\ in Math\/}.\ 
 {\bf 5}
 (1970),
349--364. 


\bibitem{OW} \name{D. Ornstein} and \name{B. Weiss}, 
Finitely determined implies very weak Bernoulli, 
{\it Israel J.\ Math\/}.\ 
{\bf 17}
(1974),
94--104.

\bibitem{OW1} \bibline,
Geodesic flows are Bernoullian, 
{\it Israel J.\ Math\/}.\ 
{\bf 14}
(1973),
184--198.


\end{references}
\end{document}